\documentclass{kluwer}    

\newdisplay{guess}{Conjecture}

\begin{document}
\begin{article}
\begin{opening}
\title{On the flag curvature of invariant Randers metrics}
\author{Hamid Reza \surname{Salimi Moghaddam}\thanks{E-mail: salimi.moghaddam@gmail.com and hr.salimi@sci.ui.ac.ir}}
\runningauthor{Hamid Reza Salimi Moghaddam} \runningtitle{On the
flag curvature of invariant Randers metrics} \institute{Department of Mathematics, Faculty of  Sciences, University of Isfahan, Isfahan,81746-73441-Iran.}

\begin{abstract}
In the present paper, the flag curvature of invariant Randers
metrics on homogeneous spaces and Lie groups is studied. We first
give an explicit formula for the flag curvature of invariant
Randers metrics arising from invariant Riemannian metrics on
homogeneous spaces and, in special case, Lie groups. We then study
Randers metrics of constant positive flag curvature and complete
underlying Riemannian metric on Lie groups. Finally we give some
properties of those Lie groups which admit a left invariant
non-Riemannian Randers metric of Berwald type arising from a left
invariant Riemannian metric and a left invariant vector field.
\end{abstract}

\classification{AMS 2000 Mathematics Subject Classification}{
22E60, 53C60, 53C30.}

\end{opening}

\section{Introduction}
The geometry of invariant structures on homogeneous spaces is one
of the interesting subjects in differential geometry. Invariant
metrics are of these invariant structures. K. Nomizu studied many
interesting properties of invariant Riemannian metrics and the
existence and properties of invariant affine connections on
reductive homogeneous spaces (see \cite{KoNo,No}.). Also some
curvature properties of invariant Riemannian metrics on Lie groups
has studied by J. Milnor \cite{Mi}. So it is important to study
invariant Finsler metrics which are a generalization of invraint
Riemannian metrics.

S. Deng and Z. Hou studied invariant Finsler metrics on reductive
homogeneous spaces and gave an algebraic description of these
metrics \cite{DeHo1,DeHo2}. Also, in \cite{EsSa1,EsSa2}, we have
studied the existence of invariant Finsler metrics on quotient
groups and the flag curvature of invariant Randers metrics on
naturally reductive homogeneous spaces. In this paper we study the
flag curvature of invariant Randers metrics on homogeneous spaces
and Lie groups. Flag curvature, which is a generalization of the
concept of sectional curvature in Riemannian geometry, is one of
the fundamental quantities which associate with a Finsler space.
In general, the computation of the flag curvature of Finsler
metrics is very difficult, therefore it is important to find an
explicit and applicable formula for the flag curvature. One of
important Finsler metrics which have found many applications in
physics are Randers metrics (see \cite{AnInMa,As}.). In this
article, by using P\"uttmann's formula \cite{Pu}, we give an
explicit formula for the flag curvature of invariant Randers
metrics arising from invariant Riemannian metrics on homogeneous
spaces and Lie groups. Then the Randers metrics of constant
positive flag curvature and complete underlying Riemannian metric
on Lie groups are studied. Finally we give some properties of
those Lie groups which admit a left invariant non-Riemannian
Randers metric of Berwald type arising from a left invariant
Riemannian metric and a left invariant vector field.

\section{Flag curvature of invariant Randers metrics on homogeneous spaces}
The aim of this section is to give an explicit formula for the
flag curvature of invariant Randers metrics of Berwald type,
arising from invariant Riemannian metrics, on homogeneous spaces.
For this purpose we need the P\"uttmann's formula for the
curvature tensor of invariant Riemannian metrics on homogeneous
spaces (see \cite{Pu}.).

Let $G$ be a compact Lie group, $H$ a closed subgroup, and $g_0$ a
bi-invariant Riemannian metric on $G$. Assume that $\frak{g}$ and
$\frak{h}$ are the Lie algebras of $G$ and $H$ respectively. The
tangent space of the homogeneous space $G/H$ is given by the
orthogonal compliment $\frak{m}$ of $\frak{h}$ in $\frak{g}$ with
respect to $g_0$. Each invariant metric $g$ on $G/H$ is determined
by its restriction to $\frak{m}$. The arising $Ad_H$-invariant
inner product from $g$ on $\frak{m}$ can extend to an
$Ad_H$-invariant inner product on $\frak{g}$ by taking $g_0$ for
the components in $\frak{h}$. In this way the invariant metric $g$
on $G/H$ determines a unique left invariant metric on $G$ that we
also denote by $g$. The values of $g_0$ and $g$ at the identity
are inner products on $\frak{g}$ which we denote as $<.,.>_0$ and
$<.,.>$. The inner product $<.,.>$ determines a positive definite
endomorphism $\phi$ of $\frak{g}$ such that $<X,Y>=<\phi X,Y>_0$
for all $X, Y\in\frak{g}$.\\

Now we give the following Lemma which was proved by T. P\"uttmann
(see \cite{Pu}.).

\newtheorem{lem}{Lemma}
\begin{lem}\label{Puttmann}
The curvature tensor of the invariant metric $<.,.>$ on the
compact homogeneous space $G/H$ is given by
\begin{eqnarray}\label{puttmans formula}
  <R(X,Y)Z,W> &=& \frac{1}{2}(<B_-(X,Y),[Z,W]>_0+<[X,Y],B_-(Z,W)>_0) \nonumber \\
    &+& \frac{1}{4}(<[X,W],[Y,Z]_{\frak{m}}>-<[X,Z],[Y,W]_{\frak{m}}> \\
    &-& 2<[X,Y],[Z,W]_{\frak{m}}>)+(<B_+(X,W),\phi^{-1}B_+(Y,Z)>_0 \nonumber\\
    &-&<B_+(X,Z),\phi^{-1}B_+(Y,W)>_0)\nonumber,
\end{eqnarray}
where the symmetric resp.skew symmetric bilinear maps $B_+$ and
$B_-$ are defined by
\begin{eqnarray*}
  B_+(X,Y) &=& \frac{1}{2}([X,\phi Y]+[Y,\phi X]), \\
  B_-(X,Y) &=& \frac{1}{2}([\phi X,Y]+[X,\phi Y]),
\end{eqnarray*}
and $[.,.]_{\frak{m}}$ is the projection of $[.,.]$ to
$\frak{m}$.\hfill$\Box$
\end{lem}

Let $\tilde{X}$ be an invariant vector field on the homogeneous
space $G/H$ such that
$\parallel\tilde{X}\parallel=\sqrt{g(\tilde{X},\tilde{X})}<1$. A
case happen when $G/H$ is reductive with
$\frak{g}=\frak{m}\oplus\frak{h}$ and $\tilde{X}$ is the
corresponding left invariant vector field to a vector
$X\in\frak{m}$ such that $<X,X><1$ and $Ad(h)X=X$ for all $h\in H$
(see \cite{DeHo2} and \cite{EsSa1}.). By using $\tilde{X}$ we can
construct an invariant Randers metric on the homogeneous space
$G/H$ in the following way:
\begin{eqnarray}
  F(xH,Y) = \sqrt{g(xH)(Y,Y)}+g(xH)(\tilde{X}_x,Y) \ \ \ \ \forall Y\in
  T_{xH}(G/H).
\end{eqnarray}
Now we give an explicit formula for the flag curvature of these
invariant Randers metrics.

\newtheorem{thm}{Theorem}
\begin{thm}\label{flagcurvature}
Let $G$ be a compact Lie group, $H$ a closed subgroup, $g_0$ a
bi-invariant metric on $G$, and $\frak{g}$ and $\frak{h}$ the Lie
algebras of $G$ and $H$ respectively. Also let $g$ be any
invariant Riemannian metric on the homogeneous space $G/H$ such
that $<Y,Z>=<\phi Y,Z>_0$ for all $Y, Z\in \frak{g}$. Assume that
$\tilde{X}$ is an invariant vector field on $G/H$ which is
parallel with respect to $g$ and $g(\tilde{X},\tilde{X})<1$ and
$\tilde{X}_H=X$. Suppose that $F$ is the Randers metric arising
from $g$ and $\tilde{X}$, and $(P,Y)$ is a flag in $T_H(G/H)$ such
that $\{Y,U\}$ is an orthonormal basis of $P$ with respect to
$<.,.>$. Then the flag curvature of the flag $(P,Y)$ in $T_H(G/H)$
is given by
\begin{eqnarray}\label{Fcurvatureformula}
  K(P,Y)=\frac{A}{(1+<X,Y>)^2(1-<X,Y>)},
\end{eqnarray}
where $A=\alpha.<X,U>+\gamma(1+<X,Y>$, and for $A$ we have:
\begin{eqnarray}
    \alpha&=&\frac{1}{4}(<[\phi U,Y]+[U,\phi Y],[Y,X]>_0+<[U,Y],[\phi Y,X]+[Y,\phi X]>_0)\nonumber\\
    &&+\frac{3}{4}<[Y,U],[Y,X]_\frak{m}>+\frac{1}{2}<[U,\phi X]+[X,\phi U],\phi^{-1}([Y,\phi Y])>_0\nonumber\\
    &&-\frac{1}{4}<[U,\phi Y]+[Y,\phi U],\phi^{-1}([Y,\phi X]+[X,\phi
    Y])>_0,
\end{eqnarray}
and
\begin{eqnarray}
  \gamma&=&\frac{1}{2}<[\phi U,Y]+[U,\phi Y],[Y,X]>_0\nonumber \\
  && \ \ \ +\frac{3}{4}<[Y,U],[Y,U]_{\frak{m}}>+<[U,\phi U],\phi^{-1}([Y,\phi Y])>_0 \\
  && \ \ \ -\frac{1}{4}<[U,\phi Y]+[Y,\phi U],\phi^{-1}([Y,\phi U]+[U, \phi Y])>_0.\nonumber
\end{eqnarray}
\end{thm}

\newproof{pf}{Proof}
\begin{pf}
$\tilde{X}$ is parallel with respect to $g$, therefore $F$ is of
Berwald type and the Chern connection of $F$ and the Riemannian
connection of $g$ coincide (see \cite{BaChSh}, page 305.), so we
have $R^F(U,V)W=R^g(U,V)W$, where $R^F$ and $R^g$ are the
curvature tensors of $F$ and $g$, respectively. Let $R:=R^g=R^F$
be the curvature tensor of $F$ (or $g$). Also for the flag
curvature we have (\cite{Sh}):
\begin{equation}\label{flag}
   K(P,Y)=\frac{g_Y(R(U,Y)Y,U)}{g_Y(Y,Y).g_Y(U,U)-g_Y^2(Y,U)},
\end{equation}
where $g_Y(U,V)=\frac{1}{2}\frac{\partial^2}{\partial s\partial
t}(F^2(Y+sU+tV))|_{s=t=0}$.\\
By a direct computation for $F$ we get
\begin{eqnarray}\label{g_Y}
    g_Y(U,V)&=&g(U,V)+g(X,U).g(X,V)-\frac{g(X,Y).g(Y,V).g(Y,U)}{g(Y,Y)^{\frac{3}{2}}}+\nonumber\\
    &&\frac{1}{\sqrt{g(Y,Y)}}\{g(X,U).g(Y,V)+g(X,Y).g(U,V)\\
    &&+g(X,V).g(Y,U)\}.\nonumber
\end{eqnarray}
Since $\{Y,U\}$ is an orthonormal basis of $P$ with respect to
$<.,.>$, by using the formula \ref{g_Y} we have:
\begin{eqnarray}\label{eq1}
  g_Y(Y,Y).g_Y(U,U)-g_Y(Y,U)=(1+<X,Y>)^2(1-<X,Y>).
\end{eqnarray}
Also we have:
\begin{eqnarray}\label{eq2}
  g_Y(R(U,Y)Y,U)&=&<R(U,Y)Y,U>+<X,R(U,Y)Y>.<X,U>\nonumber\\
  &&+<X,Y>.<R(U,Y)Y,U>\\&&
  +<X,U>.<Y,R(U,Y)Y>,\nonumber
\end{eqnarray}
now let $\alpha=<X,R(U,Y)Y>$, $\theta=<Y,R(U,Y)Y>$ and
$\gamma=<R(U,Y)Y,U>$.\\
By using P\"uttmann's formula (see Lemma \ref{Puttmann}.) and some
computations we have:
\begin{eqnarray}\label{eq3}
\alpha&=&\frac{1}{4}(<[\phi U,Y]+[U,\phi Y],[Y,X]>_0+<[U,Y],[\phi Y,X]+[Y,\phi X]>_0)\nonumber\\
    &&+\frac{3}{4}<[Y,U],[Y,X]_\frak{m}>+\frac{1}{2}<[U,\phi X]+[X,\phi U],\phi^{-1}([Y,\phi Y])>_0\nonumber\\
    &&-\frac{1}{4}<[U,\phi Y]+[Y,\phi U],\phi^{-1}([Y,\phi X]+[X,\phi
    Y])>_0,
\end{eqnarray}
\begin{eqnarray}\label{eq4}
  \theta=0,
\end{eqnarray}
and
\begin{eqnarray}\label{eq5}
      \gamma&=&\frac{1}{2}<[\phi U,Y]+[U,\phi Y],[Y,U]>_0+\frac{3}{4}<[Y,U],[Y,U]_{\frak{m}}>\nonumber\\
             &&+<[U,\phi U],\phi^{-1}([Y,\phi Y])>_0\\
             && -\frac{1}{4}<[U,\phi Y]+[Y,\phi U],\phi^{-1}([Y,\phi U]+[U, \phi Y])>_0.\nonumber
\end{eqnarray}
Substituting the equations (\ref{g_Y}), (\ref{eq1}), (\ref{eq2}),
(\ref{eq3}), (\ref{eq4}) and (\ref{eq5}) in the equation
(\ref{flag}) completes the proof. \hfill$\Box$
\end{pf}

\newproof{rem}{Remark}
\begin{rem}
In the previous theorem, If we let $H=\{e\}$ and
$\frak{m}=\frak{g}$ then we can obtain a formula for the flag
curvature of the left invariant Randers metrics of Berwald types
arising from a left invariant Riemannian metric $g$ and a left
invariant vector field $\tilde{X}$ on Lie group $G$.
\end{rem}

If the invariant Randers metric arises from a bi-invariant
Riemannian metric on a Lie group then we can obtain a simpler
formula for the flag curvature, we give this formula in the
following theorem.

\begin{thm}
Suppose that $g_0$ is a bi-invariant Riemannian metric on a Lie
group $G$ and $\tilde{X}$ is a left invariant vector field on $G$
such that $g_0(\tilde{X},\tilde{X})<1$ and $\tilde{X}$ is parallel
with respect to $g_0$. Then we can define a left invariant Randers
metric $F$ as follows:
\begin{eqnarray*}
  F(x,Y)=\sqrt{g_0(x)(Y,Y)}+g_0(x)(\tilde{X}_x,Y).
\end{eqnarray*}
Assume that $(P,Y)$ is a flag in $T_eG$ such that $\{Y,U\}$ is an
orthonormal basis of $P$ with respect to $<.,.>_0$. Then the flag
curvature of the flag $(P,Y)$ in $T_eG$ is given by
\begin{eqnarray*}
  K(P,Y)=\frac{<[Y,[U,Y]],X>_0.<X,U>_0+<[Y,[U,Y]],U>_0(1+<X,Y>_0)}{4(1+<X,Y>_0)^2(1-<X,Y>_0)}.
\end{eqnarray*}
\end{thm}

\begin{pf}
Since $\tilde{X}$ is parallel with respect to $g_0$ the curvature
tensors of $g_0$ and $F$ coincide. On the other hand for $g_0$ we
have $R(X,Y)Z=\frac{1}{4}[Z,[X,Y]]$, therefore by substituting $R$
in the equation (\ref{flag}) and using equation (\ref{g_Y}) the
proof is completed. \hfill$\Box$
\end{pf}

\section{Invariant Randers metrics on Lie groups}
In this section we study the left invariant Randers metrics on Lie
groups and, in some special cases, find some results about the
dimension of Lie groups which can admit invariant Randers metrics.
These conclusions are obtained by using Yasuda-Shimada theorem.
The Yasuda-Shimada theorem is one of important theorems which
characterize the Randers spaces. In the year 2001, Shen's examples
of Randers manifolds with constant flag curvature motivated Bao
and Robles to determine necessary and sufficient conditions for a
Randers manifold to have constant flag curvature. Shen's examples
showed that the original version of Yasuda-Shimada theorem (1977)
is wrong. Then Bao and Robles corrected the Yasuda-Shimada theorem
(1977) and gave the correct version of this theorem,
Yasuda-Shimada theorem (2001) (see \cite{BaRo}.). (For a
comprehensive history of Yasuda-Shimada theorem see \cite{Ba}.)\\

Suppose that $M$ is an $n$-dimensional manifold endowed with a
Riemannian metric $g=(g_{ij}(x))$ and a nowhere zero 1-form
$b=(b_i(x))$ such that $\|b\|=b_i(x)b_j(x)g^{ij}(x)<1$. We can
define a Randers metric on $M$ as follows
\begin{equation}\label{eq6}
    F(x,Y)=\sqrt{g_{ij}(x)Y^iY^j}+b_i(x)Y^i.
\end{equation}
Next, we consider the 1-form $\beta=b^i(b_{j|i}-b_{i|j})dx^i$,
where the covariant derivative is taken with respect to
Levi-Civita connection to $M$. Now we give the Yasuda-Shimada
theorem from \cite{Ba}.

\begin{thm}\label{Yasuda-Shimada}
(Yasuda-Shimada) (see \cite{Ba}.) Let $F$ be a strongly convex
non-Riemannian Randers metric on a smooth manifold $M$ of
dimension $n\geq 2$. Let $g_{ij}$ be the underlying Riemannian
metric and $b_i$ the drift 1-form. Then:

\begin{description}
    \item[(+)] $F$ satisfies $\beta=0$ and has constant positive
    flag curvature $K$ if and only if:
    \begin{itemize}
        \item $b$ is a non-parallel Killing field of $g$ with
        constant length;
        \item the Riemann curvature tensor of $g$ is given by
        \begin{eqnarray*}
         R_{hijk}&=&K(1-\|b\|^2)(g_{hk}g_{ij}-g_{hj}g_{ik})\\
         &&+K(g_{ij}b_hb_k-g_{ik}b_hb_j)\\
         &&-K(g_{hj}b_ib_k-g_{hk}b_ib_j)\\
         &&-b_{i|j}b_{h|k}+b_{i|k}b_{h|j}+2b_{h|i}b_{j|k}
        \end{eqnarray*}
    \end{itemize}
    \item[(0)] $F$ satisfies $\beta=0$ and has zero flag
    curvature $\Leftrightarrow$ it is locally Minkowskian.
    \item[(--)] $F$ satisfies $\beta=0$ and has constant negative
    flag curvature if and only if:
    \begin{itemize}
        \item $b$ is a closed 1-form;
        \item $b_{i|k}=\frac{1}{2}\sigma(g_{ik}-b_ib_k)$, with
        $\sigma^2=-16K$;
        \item $g$ has constant negative sectional curvature $4K$,
        that is, \\ $R_{hijk}=4K(g_{ij}g_{hk}-g_{ik}g_{hj})$.
    \end{itemize}
\end{description} \hfill$\Box$
\end{thm}

Since any Randers manifold of dimension $n=1$ is a Riemannian
manifold from now on we consider $n>1$.

An immediate conclusion of Yasuda-Shimada theorem is the following
corollary.

\newtheorem{cor}{Corollary}
\begin{cor}
There is no non-Riemannian Randers metric of Berwald type with
$\beta=0$ and constant positive flag curvature.
\end{cor}

Now by using the results of \cite{BeFa} we obtain the following
conclusions.

\begin{thm}\label{parallel}
Let $F^n=(M,F,g_{ij},b_i)$ be an $n$-dimensional parallelizable
Randers manifold of constant positive flag curvature with
$\beta=0$ on $M$ and complete Riemannian metric $g=(g_{ij})$. Then
the dimension of $M$ must be $3$ or $7$.
\end{thm}

\begin{pf}
By using theorem 2.2 of \cite{BeFa} $M$ is diffeomorphic with a
sphere of dimension $n=2k+1$. But a sphere $S^m$ is parallelizable
if and only if $m=1,3$ or $7$ (see \cite{Ad}.). Therefore $n=3$ or
$7$. \hfill$\Box$
\end{pf}

A family of Randers metrics of constant positive flag curvature on
Lie group $S^3$ was studied by D. Bao and Z. Shen (see
\cite{BaSh}.). They produced, for each $K>1$, an explicit example
of a compact boundaryless (non-Riemannian) Randers spaces that has
constant positive flag curvature $K$, and which is not
projectively flat, on Lie group $S^3$. In the following we give
some results about the dimension of Lie groups which can admit
Randers metrics of constant positive flag curvature. These results
show that the dimension $3$ is important.

\begin{cor}
There is no Randers Lie group of constant positive flag curvature
with $\beta=0$, complete Riemannian metric $g=(g_{ij})$ and $n\neq
3$.
\end{cor}

\begin{pf}
Any Lie group is parallelizable, so by attention to theorem
\ref{parallel} and the condition $n\neq 3$, $n$ must be $7$. Since
$G$ is diffeomorphic to $S^7$ and $S^7$ can not admit any Lie
group structure, hence the proof is completed. \hfill$\Box$
\end{pf}

Similar to the \cite{Mi} for the sectional curvature of the left
invariant Riemannian metrics on Lie groups, we compute the flag
curvature of the left invariant Randers metrics on Lie groups in
the following theorem.

\begin{thm}\label{flag(ei)}
Let $G$ be a compact Lie group with Lie algebra $\frak{g}$, $g_0$
a bi-invariant Riemannian metric on $G$, and $g$ any left
invariant Riemannian metric on $G$ such that $<X,Y>=<\phi X,Y>_0$
for a positive definite endomorphism
$\phi:\frak{g}\longrightarrow\frak{g}$. Assume that $X\in\frak{g}$
is a vector such that $<X,X><1$ and $F$ is the Randers metric
arising from $\tilde{X}$ and $g$ as follows:
\begin{eqnarray*}
  F(x,Y)=\sqrt{g(x)(Y,Y)}+g(x)(\tilde{X}_x,Y),
\end{eqnarray*}

where $\tilde{X}$ is the left invariant vector field corresponding
to $X$, and we have assumed $\tilde{X}$ is parallel with respect
to $g$. Let $\{e_1,\cdots,e_n\}\subset\frak{g}$ be a
$g$-orthonormal basis for $\frak{g}$. Then the flag curvature of
$F$ for the flag $P=span\{e_i,e_j\} (i\neq j)$ at the point
$(e,e_i)$, where $e$ is the unit element of $G$, is given by the
following formula:
\begin{eqnarray*}
 K(P=span\{e_i,e_j\},e_i)=\frac{X_j.<R(e_j,e_i)e_i,X>+(1+X_i).<R(e_j,e_i)e_i,e_j>}{(1+X_i)^2(1-X_i)},
\end{eqnarray*}
 where $X=X^ke_k$,
 \begin{eqnarray*}
   <R(e_j,e_i)e_i,X>&=& -\frac{1}{4}(<[\phi e_j,e_i],[e_i,X]>_0+<[e_j,\phi e_i],[e_i,X]>_0 \\
   &&+<[e_j,e_i],[\phi e_i,X]>_0+<[e_j,e_i],[e_i,\phi X]>_0)\\
   &&+\frac{3}{4}<[e_j,e_i],[e_i,X]>\\
   &&-\frac{1}{2}<[e_j,\phi X]+[X,\phi e_j],\phi^{-1}([e_i,\phi e_i])>_0\\
   &&+\frac{1}{4}<[e_j,\phi e_i]+[e_i,\phi e_j],\phi^{-1}([e_i,\phi X]+[X,\phi e_i])>_0
 \end{eqnarray*}
 and
 \begin{eqnarray*}
   <R(e_j,e_i)e_i,e_j>&=&-\frac{1}{2}(<[\phi e_j,e_i],[e_i,e_j]>_0+<[e_j,\phi e_i],[e_i,e_j]>_0) \\
    &&+\frac{3}{4}<[e_j,e_i],[e_i,e_j]>-<[e_j,\phi e_j],\phi^{-1}([e_i,\phi e_i])>_0\\
    &&+\frac{1}{4}<[e_j,\phi e_i]+[e_i,\phi e_j],\phi^{-1}([e_i,\phi e_j]+[e_j,\phi e_i])>_0.
 \end{eqnarray*}
\end{thm}

\begin{pf}
By using theorem \ref{flagcurvature}, the proof is clear.
\hfill$\Box$
\end{pf}

Now we give some properties of those Lie groups which admit a left
invariant non-Riemannian Randers metric of Berwald type arising
from a left invariant Riemannian metric and a left invariant
vector field.

\begin{thm}
There is no left invariant non-Riemannian Randers metric of
Berwald type arising from a left invariant Riemannian metric and a
left invariant vector field on connected Lie groups with a perfect
Lie algebra, that is, a Lie algebra $\frak{g}$ for which the
equation $[\frak{g},\frak{g}]=\frak{g}$ holds.
\end{thm}

\begin{pf}
If a left invariant vector field $X$ is parallel with respect to a
left invariant Riemannian metric $g$ then, by using Lemma 4.3 of
\cite{BrFiSpTaWu}, $g(X,[\frak{g},\frak{g}])=0$. Since $\frak{g}$
is perfect therefore $X$ must be zero. \hfill$\Box$
\end{pf}

\begin{cor}
There is not any left invariant non-Riemannian Randers metric of
Berwald type arising from a left invariant Riemannian metric and a
left invariant vector field on semisimple connected Lie groups.
\end{cor}

\begin{cor}
If a Lie group $G$ admits a left invariant non-Riemannian Randers
metric of Berwald type $F$ arising from a left invariant
Riemannian metric $g$ and a left invariant vector field $X$ then
for sectional curvature of the Riemannian metric $g$ we have
\begin{eqnarray*}
  K(X,u)\geqslant 0
\end{eqnarray*}
for all $u$, where equality holds if and only if $u$ is orthogonal
to the image $[X,\frak{g}]$.
\end{cor}

\begin{pf}
Since $F$ is of Berwald type, $X$ is parallel with respect to $g$.
By using Lemma 4.3 of \cite{BrFiSpTaWu}, $ad(X)$ is skew-adjoint,
therefore by Lemma 1.2 of \cite{Mi} we have $K(X,u)\geqslant 0$.
\hfill$\Box$
\end{pf}


\end{article}
\end{document}